\documentclass{book}    

\usepackage[numbers]{natbib}

\usepackage{piers}  
\begin{document}

\title{Variational Integrators for Maxwell's Equations with Sources}
\maketitle

\author      {F. M. Lastname}
\affiliation {University}
\address     {}
\city        {Boston}
\postalcode  {}
\country     {USA}
\phone       {345566}    
\fax         {233445}    
\email       {email@email.com}  
\misc        { }  
\nomakeauthor

\author      {F. M. Lastname}
\affiliation {University}
\address     {}
\city        {Boston}
\postalcode  {}
\country     {USA}
\phone       {345566}    
\fax         {233445}    
\email       {email@email.com}  
\misc        { }  
\nomakeauthor

\begin{authors}

{\bf A. Stern}$^{1}$, {\bf Y. Tong}$^{1,2}$, {\bf M. Desbrun}$^{1}$, {\bf and J. E. Marsden}$^{1}$\\
\medskip
$^{1}$California Institute of Technology, USA\\
$^{2}$Michigan State University, USA

\end{authors}

\begin{paper}

\begin{piersabstract}
  In recent years, two important techniques for geometric numerical
  discretization have been developed.  In computational
  electromagnetics, spatial discretization has been improved by the
  use of mixed finite elements and discrete differential forms.
  Simultaneously, the dynamical systems and mechanics communities have
  developed structure-preserving time integrators, notably variational
  integrators that are constructed from a Lagrangian action principle.
  Here, we discuss how to combine these two frameworks to develop {\em
    variational spacetime integrators} for Maxwell's equations.
  Extending our previous work, which first introduced this variational
  perspective for Maxwell's equations without sources, we also show
  here how to incorporate free sources of charge and current.
\end{piersabstract}

\psection{Introduction}

In computational electromagnetics, as in an increasing number of other
fields in applied science and engineering, there is both practical and
theoretical interest in developing {\em geometric numerical
  integrators}.  These numerical methods preserve, by construction,
various geometric properties and invariants of the continuous physical
systems that they approximate.  This is particularly important for
applications where even high-order methods may fail to capture
important features of the underlying dynamics. In this short paper, we
show that the traditional Yee scheme and extensions can be derived
from the Euler-Lagrange equations of a discrete action, i.e., by
designing an electromagnetic \emph{variational integrator}, including
free sources of charge and current in non-dissipative
media. Furthermore, we present how to use this discrete geometric
framework to allow for asynchronous time stepping on unstructured
grids, as recently introduced in~\citet{StToDeMa2007}.

Variational integrators (not to be confused with variational methods such as
finite element schemes) were originally developed for geometric time
integration, particularly to simulate dynamical systems in Lagrangian
mechanics. The key idea is the following: rather than approximating the
equations of motion directly, one discretizes the Lagrangian and its
associated action integral (e.g., using a numerical quadrature rule), and then
{\em derives} a structure-preserving approximation to the equations of motion
by applying Hamilton's principle of stationary action. Since the numerical
method is derived from a Lagrangian variational principle, some important
results from Lagrangian dynamics carry over to the discretized system,
including Noether's theorem relating symmetries to conserved momentum maps, as
well as the fact that the Euler-Lagrange flow is a symplectic
mapping. (\citet[See][]{MaWe2001,LeMaOrWe2004}.)

\paragraph{Overview.}
To develop a variational integrator for Maxwell's equations, the
discrete Hamilton's principle needs to incorporate more than just the
time discretization, as in mechanics; spatial discretization also
needs to be handled carefully. Building upon mixed finite elements in
space~\citep{Nedelec1980,Bossavit1998,GrKo2004}, we treat the
electromagnetic {\em Lagrangian density} as a discrete differential
$4$-form in spacetime.  Extremizing the integral of this Lagrangian
density with fixed boundary conditions directly leads to discrete
update rules for the electromagnetic fields, with either uniform or
asynchronous time steps across the various spatial elements.

\psection{Review of Maxwell's Equations in Spacetime}

\paragraph{Electromagnetic Forms.}  Let $A$ be a $1$-form on spacetime,
called the {\em electromagnetic potential}, and then define the {\em
  Faraday $2$-form} to be its exterior derivative $ F = d A $.  Given
a time coordinate $t$, this splits into the components
\begin{equation*}
F = E \wedge d t + B ,
\end{equation*}
where $E$ is the electric displacement $1$-form and $B$ is the
magnetic flux $2$-form, both defined on the spacelike Cauchy surfaces
$\Sigma$ with constant $t$.  If $ * $ is the Hodge star associated to
the spacetime metric, then we can also split the dual $2$-form
\begin{equation*}
  *F = \left( *_\mu B \right) \wedge d t - * _\epsilon E = H \wedge d
  t - D ,
\end{equation*}
where (again, restricted to Cauchy surfaces) $H$ is the magnetic
displacement $1$-form, $D$ is the electric flux $2$-form, and $ * _\mu
$ and $ * _\epsilon $ are respectively the magnetic permeability and
electric permittivity.  Finally, for systems with free sources, there
is a {\em source 3-form} $\mathcal{J}$, satisfying the continuity of
charge condition $ d \mathcal{J} = 0 $. In terms of coordinates, this
can be split into
\begin{equation*}
\mathcal{J} = J \wedge d t - \rho ,
\end{equation*}
where $J$ is the current density $2$-form and $\rho$ is the charge
density $3$-form on Cauchy surfaces.

\paragraph{Maxwell's Equations.}  With the spacetime forms and
operators defined above, Maxwell's equations become
\begin{equation*} 
  d F = 0, \qquad 
  d {*F} = \mathcal{J} .
\end{equation*} 
Note that the first equation follows automatically from $ F = d A $, since
taking the exterior derivative of both sides yields $ d F = dd A = 0 $. The
second equation is consistent with the continuity of charge condition, since $
d \mathcal{J} = d d {*F} = 0 $.

\paragraph{Lagrangian Formulation.}  Given the electromagnetic
potential $1$-form $A$ and source $3$-form $\mathcal{J}$, we can define the
Lagrangian density to be the $4$-form
\begin{equation*}
\mathcal{L} = - \frac{1}{2} d A \wedge * d A + A \wedge \mathcal{J} ,
\end{equation*}
with the associated action functional $ S[A] = \int _X \mathcal{L} $
taken over the spacetime domain $X$. Suppose that $\alpha$ is a
variation of $A$, vanishing on the boundary $\partial X $. Varying the
action along $\alpha$ yields
\begin{equation*} 
\mathbf{d} S[A] \cdot \alpha = \int _X \left( - d \alpha \wedge *d A +
  \alpha \wedge \mathcal{J} \right)
= \int _X \alpha \wedge \left( -d {*d} A + \mathcal{J} \right) .
\end{equation*} 
Hamilton's principle of stationary action states that this variation must
equal zero for any such $\alpha$, implying the Euler-Lagrange equations
$d{*d}A = \mathcal{J} $. Finally, substituting $ F = d A $ and recalling that
$ d F = d d A = 0 $, we see that this is equivalent to Maxwell's equations.

\psection{Geometric Properties of Maxwell's Equations}

As written in terms of $F$ above, Maxwell's equations have 8
components: 6 dynamical equations, which describe how the fields
change in time, and 2 ``divergence constraints'' containing only
spatial derivatives.  The fact that these constraints are
automatically preserved by the dynamical equations (and can therefore
effectively be ignored except at the initial time) comes directly from
the differential gauge symmetry and Lagrangian variational structure.
We discuss these geometric properties here, with a view towards
developing numerical methods that preserve them.

\paragraph{Reduction by Gauge Fixing.}  Maxwell's equations are
invariant under gauge transformations $ A \mapsto A + d f $ for any
scalar function $f$, since taking the exterior derivative maps $ F
\mapsto F + d d f = F $.  Therefore, given a time coordinate $t$, we
can fix the gauge so that $ A \cdot \frac{\partial}{\partial t} = 0 $,
i.e., $A$ has only spacelike components.  This partial gauge fixing is
known as the {\em Weyl gauge}.  Restricted to this subspace of
potentials, the Lagrangian then becomes
\begin{align*}
  \mathcal{L} &= -\frac{1}{2} \left( d _t A + d _\Sigma A \right)
  \wedge
  *\left( d _t A + d _\Sigma A \right)  + A \wedge \mathcal{J} \\
  &= -\frac{1}{2} \left( d _t A \wedge * d _t A + d _\Sigma A \wedge *
    d _\Sigma A \right) + A \wedge J \wedge d t
\end{align*}
Here, we have adopted the notation $ d _t $ and $ d _\Sigma $ for the
exterior derivative taken only in time and in space, respectively; in
particular, we then have $ d _t A = E \wedge d t $ and $ d _\Sigma A =
B $.

Next, varying the action along a restricted variation $\alpha$ that
vanishes on $ \partial X $,
\begin{align}
  \mathbf{d} S[A] \cdot \alpha &= \int _X \left( d _t \alpha \wedge D
    - d _\Sigma \alpha \wedge H \wedge d t + \alpha \wedge J \wedge
    dt
  \right) \label{eqn:weylaction}\\
  &= \int _X \alpha \wedge \left( d _t D - d _\Sigma H \wedge d t + J
    \wedge d t \right) .\nonumber
\end{align}
Setting this equal to zero by Hamilton's principle, one immediately
gets Amp\`ere's law as the sole Euler-Lagrange equation.  The
divergence constraint $ d _\Sigma {D} = \rho $, corresponding to
Gauss' law, has been eliminated via the restriction to the Weyl gauge.

\paragraph{Noether's Theorem Automatically Preserves Gauss' Law.}
There are two ways that one can see why Gauss' law is automatically
preserved, even though it has been eliminated from the Euler-Lagrange
equations.  The first is to take the ``divergence'' $ d _\Sigma $ of
Amp\`ere's law, obtaining
\begin{equation*}
  0 = d _\Sigma d _t D - d _\Sigma d _\Sigma H \wedge d t + d _\Sigma J
  \wedge d t
  = d _t \left( d _\Sigma D - \rho \right) .
\end{equation*}
Therefore, if this condition holds at the initial time, then it holds
for all time.

A more ``geometric'' way to obtain this result is to use Noether's
theorem, with respect to the remaining gauge symmetry $ A \mapsto A +
d _\Sigma f $ for scalar functions $f$ on $\Sigma$.  To derive this,
let us restrict $A$ to be an Euler-Lagrange solution in the Weyl
gauge, but remove the previous requirement that variations $\alpha$ be
fixed at the initial time $ t _0 $ and final time $ t _f $.  Then,
varying the action along this new $\alpha$, the Euler-Lagrange term
disappears, but we now pick up an additional boundary term due to
integration by parts
\begin{equation*}
  \mathbf{d} S [A] \cdot \alpha = \left. \int _\Sigma \alpha \wedge
    {D} \,\right| _{t _0 } ^{ t _f }.
\end{equation*}
If we vary along a gauge transformation $ \alpha = d _\Sigma f $, then
this becomes
\begin{equation*}
  \mathbf{d} S [A] \cdot d _\Sigma f = \left. \int _\Sigma d _\Sigma
    f \wedge {D} \,\right| _{t _0 } ^{ t _f }
  = - \left. \int _\Sigma f \wedge d _\Sigma {D} \,\right| _{t _0 }
  ^{ t _f }
\end{equation*}
Alternatively, plugging $ \alpha = d _\Sigma f $
into~\autoref{eqn:weylaction}, we get
\begin{equation*}
  \mathbf{d} S [A] \cdot d _\Sigma f = \int _X d _\Sigma
  f \wedge J \wedge d t
  = - \int _X f \wedge d _\Sigma J \wedge d t
  = - \int _X f \wedge d _t \rho
  = - \left. \int _\Sigma f \wedge \rho \,\right| _{ t _0 } ^{ t _f }
  .
\end{equation*}
Since these two expressions are equal, and $f$ is an arbitrary
function, it follows that
\begin{equation*}
\left. \left( d _\Sigma D - \rho \right) \right| _{t_0}^{t_f} = 0.
\end{equation*}
This indicates that $ d _\Sigma D - \rho $ is a conserved quantity, a
momentum map, so if Gauss' law holds at the initial time, then it
holds for all subsequent times as well.

\psection{Geometric Discretization of Maxwell's Equations}

Discretizing Maxwell's equations, while preserving the geometric
properties mentioned above, can be achieved using \emph{cochains} as
discrete substitutes for differential forms, as previously done in,
e.g.,~\citet{Bossavit1998}.  Therefore, to compute Maxwell's equations,
we begin by discretizing the $2$-form $F$ on a spacetime mesh $K$: $F$
assigns a real value to each oriented $2$-face of the mesh. The
exterior derivative $d$ is discretized by the coboundary operator, so
the equation $ d F = 0 $ states that $ \left\langle d F , \sigma ^3
\right\rangle = \left\langle F, \partial \sigma ^3 \right\rangle = 0
$, where $ \sigma ^3 $ is any oriented $3$-cell in $K$ and $ \partial
\sigma ^3 $ is its $2$-chain boundary.

Next, given a discrete Hodge star
operator~\cite{BoKe1999,TaKeBo1999,Hiptmair2001,AuKu2006}, $ {*F} $ is
a $2$-form on the dual mesh $ *K $, while $\mathcal{J}$ is defined as
a discrete dual $3$-form. Then, for every dual $3$-cochain $ * \sigma
^1 $ (where $ \sigma ^1 $ is the corresponding primal edge), the
equation $ d {*F} = \mathcal{J} $ becomes
\begin{equation*}
\left\langle d {* F} , * \sigma ^1 \right\rangle = \left\langle
  *F, \partial {* \sigma ^1} \right\rangle = \left\langle \mathcal{J}
  , * \sigma ^1 \right\rangle .
\end{equation*}

When the cells $ \sigma ^3 $ and $ *\sigma ^1 $ are spacetimelike,
then these correspond to the dynamical components of Maxwell's
equations, and can be used to compute subsequent values of $F$.  When
the cells are purely spacelike, they correspond to the divergence
constraint equations. The exact expression and update of these fields
in time now depends on which type of mesh and time stepping method is
desired, as described next.

\paragraph{Uniform Time Stepping.} For uniform rectangular meshes
aligned with the $ (x,y,z,t) $ axes, we can emulate the smooth coordinate
expression of $F$ as
\begin{align*}
  F &= E _x \,dx \wedge d t + E _y \,dy \wedge d t + E _z \,dz
  \wedge d t  \\
& \quad + B _x \,dy \wedge d z + B _y \,dz \wedge d x + B _z \,dx
  \wedge d y.
\end{align*}
This suggests the following discretization of $F$, shown
in~\autoref{fig:setupEB}: store $ E _x \Delta x \Delta t $ on the $ x
t $-faces, $ E _y \Delta y \Delta t $ on the $ y t $-faces, and $ E _z
\Delta z \Delta t $ on the $ z t $-faces; likewise for $B$.  To
discretize $\mathcal{J}$, we can similarly store $ J _x \Delta y
\Delta z \Delta t $, $ J _y \Delta z \Delta x \Delta t $, $ J _z
\Delta x \Delta y \Delta t $, and $ \rho \Delta x \Delta y \Delta z $
on the corresponding dual $3$-cells.
\begin{figure}[ht]
\centerline{\includegraphics[width=0.9\linewidth]{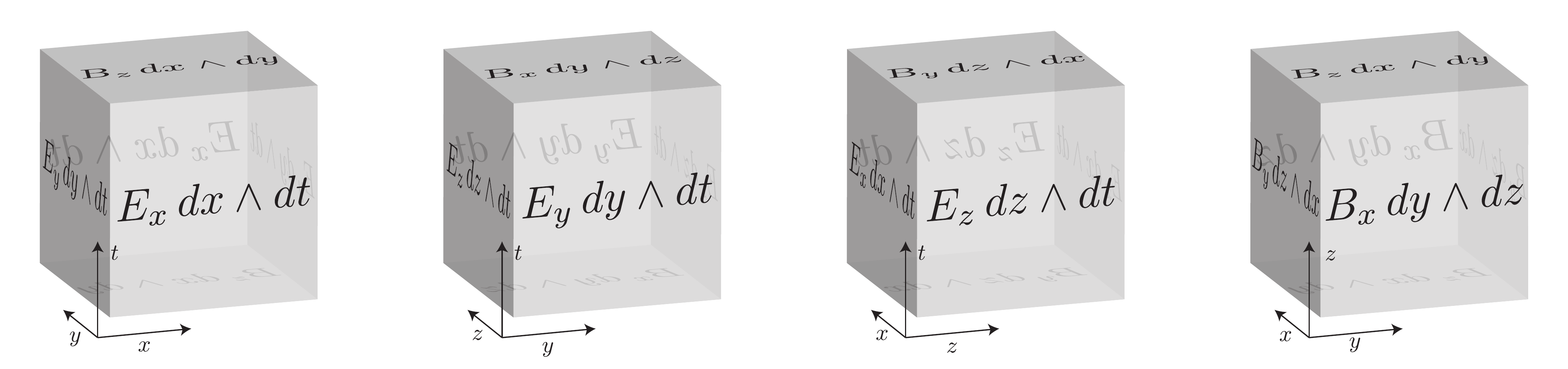}}
\caption{The $2$-form $F = E \wedge d t + B $ can be discretized, on a
  rectangular spacetime mesh, by storing the components of $E$ and $B$
  on 2D faces.  The resulting numerical method is Yee's FDTD scheme.}
    \label{fig:setupEB}
\end{figure}
If we then enforce the equations $ d F = 0 $ and $ d {*F} =
\mathcal{J} $ in the discrete sense, the result is precisely the
finite-difference time-domain (FDTD) integration scheme
of~\citet{Yee1966}. A similar procedure can be applied on unstructured
(e.g., simplicial) spatial grids, on which we take uniform time steps
$ \Delta t $ (creating prism-shaped spacetime primal elements); in
this case, solving the discrete Maxwell's equations recovers the more
recent ``Yee-like'' method of~\citet{BoKe2000}.

\paragraph{Asynchronous Time Stepping.}
As initially mentioned in~\citet{StToDeMa2007}, a more flexible
integration scheme can be designed by assigning a different time step
per element, in order to focus computational power where needed.  A
visualization of how to store $F$ on such a mesh structure is shown
in~\autoref{fig:setupAVI}.  This defines an {\em asynchronous
  variational integrator} that preserves the numerical properties of
the uniform-stepping methods outlined above. The procedure to
repeatedly update $E$ and $B$ asynchronously in time is as follows:
\begin{enumerate}
\item Select the face for which $B$ needs to be updated next.
\item $E$ advances $B$, using $ d F = 0 $.
\item $B$ advances $E$ on neighboring edges, using $ d {* F} =
  \mathcal{J} $.
\end{enumerate}
\begin{figure}[ht]
\centerline{\includegraphics[width=0.7\linewidth]{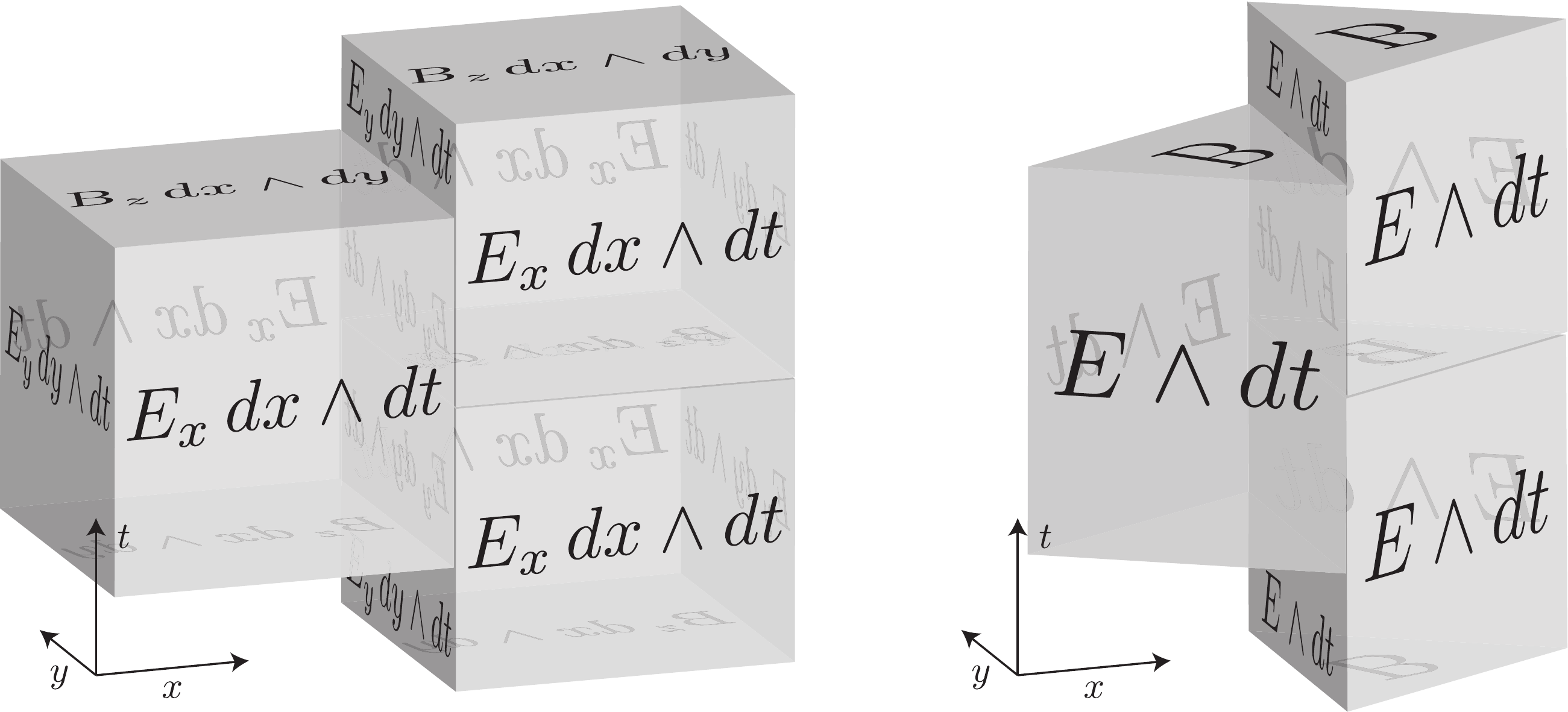}}
\caption{Here, $F$ is discretized on an asynchronous time-stepping
  grid, where each spatial element takes a different-sized time step
  from its neighbors.  This can be done for either a rectangular
  spatial grid (left) or an unstructured/simplicial spatial mesh
  (right). \label{fig:setupAVI}}
\end{figure}
Details of this algorithm, along with initial numerical results, can
be found in~\citet{StToDeMa2007}, where only the case $\mathcal{J}=0$
was described.

\begin{figure}[ht]
  \centering
  \subfigure[regular grid, uniform time steps]
  {
    \includegraphics[width=0.45\linewidth]{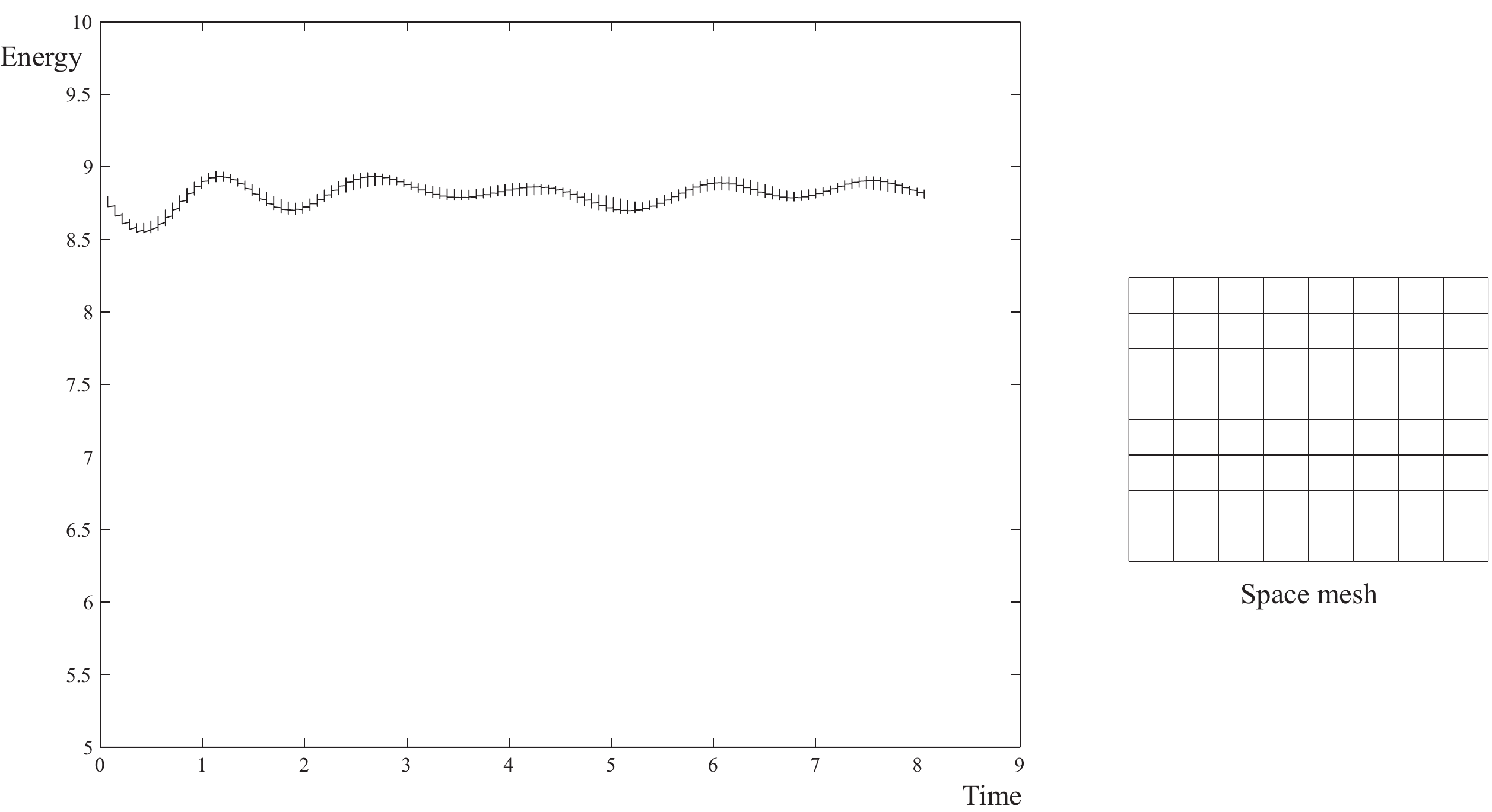}
    \label{fig:avi-energy-uniform-steps}
  }
  \subfigure[random grid, asynchronous time steps]
  {
    \includegraphics[width=0.45\linewidth]{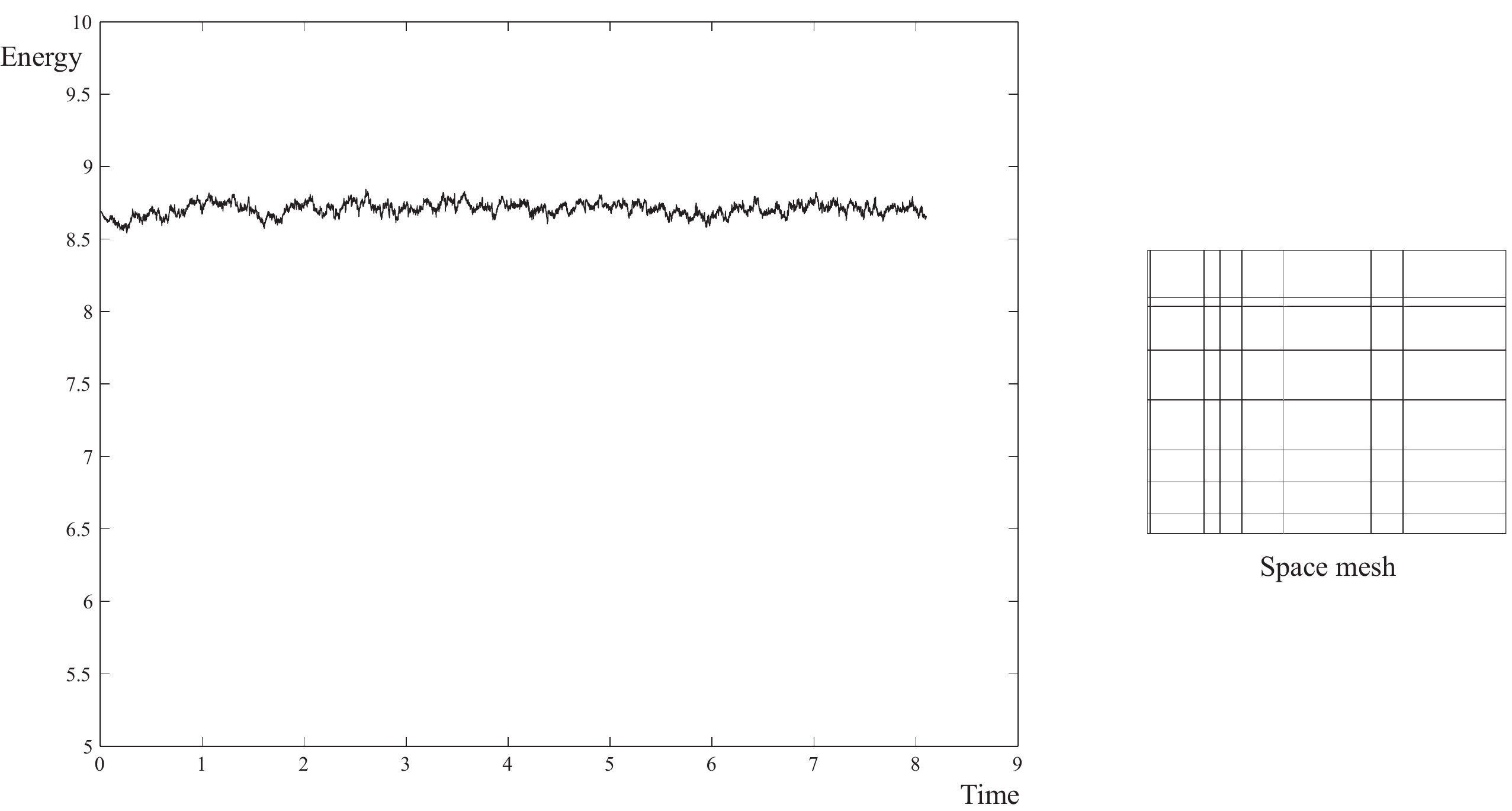}
    \label{fig:avi-energy-random-steps}
  }
  \caption{Our geometric integrator robustly maintains
    near-conservation of energy, even for asynchronous time stepping
    on a random spatial grid.}
  \label{fig:avi-energy}
\end{figure}

One promising result from the initial numerical experiments, shown
in~\autoref{fig:avi-energy}, is that this asynchronous integrator
matches the FDTD scheme's excellent energy conservation behavior, even
on a highly irregular grid, without exhibiting artificial damping or
forcing.  Future work is expected to include formal accuracy and
stability analysis of this asynchronous integrator, focusing both on
the choice of Hodge star and on resonance stability criteria for
selecting the individual time steps.

\ack Our research was partially supported by a Betty and Gordon Moore
fellowship at Caltech, NSF grants CCR-0133983 and DMS-0453145 and DOE
contract DE-FG02-04ER25657, and by NSF grant CCF-0528101.  We
gratefully acknowledge these sponsors for their support of this work.

\end{paper}


\begin{thebibliography}{28}
\providecommand{\natexlab}[1]{#1}
\providecommand{\url}[1]{#1}
\providecommand{\urlprefix}{URL }
\expandafter\ifx\csname urlstyle\endcsname\relax
  \providecommand{\doi}[1]{DOI~\discretionary{}{}{}#1}\else
  \providecommand{\doi}{DOI~\discretionary{}{}{}\begingroup
  \urlstyle{rm}\Url}\fi

\bibitem[{Auchmann and Kurz(2006)}]{AuKu2006}
Auchmann, B., Kurz, S.: A geometrically defined discrete {H}odge operator on
  simplicial cells.
\newblock IEEE Trans. Magn. \textbf{42}(4), 643--646 (2006)

\bibitem[{Bossavit(1998)}]{Bossavit1998}
Bossavit, A.: Computational electromagnetism.
\newblock Electromagnetism. Academic Press Inc., San Diego, CA (1998).
\newblock Variational formulations, complementarity, edge elements

\bibitem[{Bossavit and Kettunen(1999)}]{BoKe1999}
Bossavit, A., Kettunen, L.: Yee-like schemes on a tetrahedral mesh, with
  diagonal lumping.
\newblock Int. J. Numer. Modell. \textbf{12}(1--2), 129--142 (1999)

\bibitem[{Bossavit and Kettunen(2000)}]{BoKe2000}
Bossavit, A., Kettunen, L.: Yee-like schemes on staggered cellular grids: A
  synthesis between {FIT} and {FEM} approaches.
\newblock IEEE Trans. Magn. \textbf{36}(4), 861--867 (2000)

\bibitem[{Gross and Kotiuga(2004)}]{GrKo2004}
Gross, P.W., Kotiuga, P.R.: Electromagnetic theory and computation: a
  topological approach, \emph{Mathematical Sciences Research Institute
  Publications}, vol.~48.
\newblock Cambridge University Press, Cambridge (2004)

\bibitem[{Hiptmair(2001)}]{Hiptmair2001}
Hiptmair, R.: Discrete {H}odge operators.
\newblock Numer. Math. \textbf{90}(2), 265--289 (2001)

\bibitem[{Lew et~al.(2004)Lew, Marsden, Ortiz, and West}]{LeMaOrWe2004}
Lew, A., Marsden, J.E., Ortiz, M., West, M.: Variational time integrators.
\newblock Internat. J. Numer. Methods Engrg. \textbf{60}(1), 153--212 (2004)

\bibitem[{Marsden and West(2001)}]{MaWe2001}
Marsden, J.E., West, M.: Discrete mechanics and variational integrators.
\newblock Acta Numer. \textbf{10}, 357--514 (2001)

\bibitem[{N{\'e}d{\'e}lec(1980)}]{Nedelec1980}
N{\'e}d{\'e}lec, J.C.: Mixed finite elements in $\mathbb{R}^3$.
\newblock Numer. Math. \textbf{35}(3), 315--341 (1980)

\bibitem[Stern et~al.(2007)Stern, Tong, Desbrun, and Marsden]{StToDeMa2007}
Stern, A., Y.~Tong, M.~Desbrun, and J.~E. Marsden [2007], {\em Computational
  Electromagnetism with Variational Integrators and Discrete Differential
  Forms}, preprint at {\tt arXiv:0707.4470 [math.NA]}, 2007.

\bibitem[{Tarhasaari et~al.(1999)Tarhasaari, Kettunen, and
  Bossavit}]{TaKeBo1999}
Tarhasaari, T., Kettunen, L., Bossavit, A.: Some realizations of a discrete
  {H}odge operator: a reinterpretation of finite element techniques.
\newblock IEEE Trans. Magn. \textbf{35}(3), 1494--1497 (1999)

\bibitem[{Yee(1966)}]{Yee1966}
Yee, K.S.: Numerical solution of inital boundary value problems involving
  {M}axwell's equations in isotropic media.
\newblock IEEE Trans. Ant. Prop. \textbf{14}(3), 302--307 (1966)

\end{thebibliography}
\end{document}